\documentclass{article}
\usepackage[OT2,T1]{fontenc}
\DeclareSymbolFont{cyrletters}{OT2}{wncyr}{m}{n}
\DeclareMathSymbol{\Sha}{\mathalpha}{cyrletters}{"58}

\usepackage{amsbsy,amssymb,amscd,amsfonts,latexsym,amstext,delarray,
amsmath, diagrams}  \headsep=15pt
\usepackage[all,knot]{xy}
\def\Un{ \mbox{Un} }
\def\hpi{\hat{\pi}}

\def\cX{{\cal X} }

\def\bX{ \bar{X}}

\def\Spec{{ \mbox{Spec} }}

\def\Hom{{ \mbox{Hom} }}

\def\Om{{ \Omega }}

\def\O{{ {\cal O} }}

\def\ra{{ \rightarrow }}

\def\a{{ \alpha }}
\def\b{{ \beta} }
\def\g{{ \gamma }}
\def\d{{ \delta }}

\def\Un{ \mbox{Un} }

\def\hra{{ \hookrightarrow }}

\def\bs{ \backslash}

\def\Gal{{ \mbox{Gal} }}

\def\bQ{\bar{\Q}}

\def\cE{ {\cal E}}

\def\Z{{ \mathbb{Z}}}

\def\bq{\begin{quote}}
\def\eq{\end{quote}}

\def\Isom{ \mbox{Isom} }

\newtheorem{thm}{Theorem}[section]
\newtheorem{cor}[thm]{Corollary}
\newtheorem{lem}[thm]{Lemma}
\newtheorem{prop}[thm]{Proposition}

\def\Q{\mathbb{Q}}

\def\invlim{\varprojlim}
\def\dirlim{\varinjlim}

\def\be{\begin{equation}}
\def\ee{\end{equation}}

\def\D{ \Delta}

\def\bE{ \bar{E}}

\def\loc{\mbox{loc}}

\def\U{ U^{et}}

\def\b{ \bar}

\def\L{\Lambda}
\def\Exp{\mbox{Exp}}

\def\a{\alpha}

\def\l{ \lambda}

\def\bd{\begin{diagram} }
\def\ed{ \end{diagram}}

\def\b{\beta}

\def\bs{\begin{slide}}
\def\es{\end{slide} }

\def\Gmp{\mathbb{G}_{m,\Q_p}}

\def\bd{ \begin{diagram}}
\def\ed{\end{diagram} }

\def\Gm{\mathbb{G}_m}
\def\bGmp{\overline{\mathbb{G}}_{m,\Q_p}}

\def\Pux{\pi_1^{u,\Q_p}(\bX;b,x)}
\def\Pdrx{\pi_1^{DR}(X_{\Q_p};b,x)}
\def\piu{\pi_1^{u,\Q_p}(\bX,b)}
\def\mpiu{\pi_1^{u,\Q_p}(\bX,-b)}

\def\pih{\hat{\pi}_1(\bX,b)}

\def\ord{\mbox{ord}}

\def\A{\mathbb{A}}
\def\bsl{\backslash}
\def\pibmb{\hat{\pi}_1(\bX;b,-b) }
\def\bT{\bar{T}}
\def\cR{{\cal R}}
\title{Massey products for elliptic curves of rank 1 }
\author{Minhyong Kim}
\begin{document}
\maketitle
\begin{abstract}
For an elliptic curve over $\Q$ of
analytic rank 1, we use the level-two Selmer variety and secondary cohomology products to find explicit analytic defining
equations for global integral points inside the set of $p$-adic points.
\end{abstract}
The author must begin with an apology for writing on a topic so specific, so elementary,
and so well-understood as the study of elliptic curves of rank 1. Nevertheless,
it is hoped that a contribution not entirely without value or novelty
 is to be found within the theory of {\em Selmer varieties} for hyperbolic
 curves, applied to the complement $X=E\setminus \{e\}$ of the origin
 inside an elliptic curve $E$ over $\Q$ with $\ord_{s=1}L( E,s)=1$.
 That is, let $\cE$ be a regular minimal $\Z$-model for $E$ and $\cX$ be the complement in $\cE$
 of the closure of $e$.
 We choose an odd prime $p$ of good reduction for $\cE$.
 The main goal of the present inquiry  is
 \bq
 {\em to find explicit analytic equations defining $\cX(\Z)$ inside
 $\cX(\Z_p)$}.
 \eq
 The approach of this paper makes use of a {\em rigidified Massey product}
 in Galois cohomology\footnote{The reader is invited to consult \cite{deninger} for the corresponding construction
 in Deligne cohomology and \cite{sharifi} for a $\Gm$-analogue.}.
 That is    the \'etale local unipotent Albanese map
 $$\cX(\Z_p)\stackrel{j^{et}_{2,loc}}{\ra}H^1_f(G_p,U_2)$$
 to the level-two local Selmer variety (recalled below) associates to point
 $z$ a non-abelian cocycle $a(z)$, which can be broken canonically
 into two components $a(z)=a_1(z)+a_2(z)$, with
 $a_1(z)$ taking values in $U_1\simeq T_p(E)\otimes \Q_p$ and
 $a_2(z)$ in $U^3\bsl U^2\simeq \Q_p(1)$.
 Denoting by $c^p:G_p\ra \Q_p$ the logarithm of the $p$-adic cyclotomic character,
 we use a Massey triple product
 $$z\mapsto (c^p, a_1(z), a_1(z))\in H^2(G_p, \Q_p(1))\simeq \Q_p$$
 to construct a function on the local points of $\cX$.
 Recall that Massey products are secondary cohomology products
 arising in connection with associative differential graded algebras
 $(A,d)$: If we are given classes $[\a],[\b], [\g]\in H^1(A)$ such
 that
 $$[\a][\b]=0=[\b][\g],$$
 we can solve the equations
 $$dx=\a \b, \ \ \ dy=\b \g$$
 for $x,y \in A^1$.
 The element
 $$x\g+\a y \in A^2$$
 satisfies
 $$d(x\g+\a y )=dx\g -\a dy=0,$$
  so that we  obtain a class
 $$[x\g +\a y]\in H^2(A).$$
 Of course it depends on the choice of $x$ and $y$, so that the product as
 a function of the initial triple takes values in
 $$H^2(A)/[H^1(A)[\g]+[\a]H^1(A)].$$
 The elements $x$ and $y$ constitute a {\em defining system} for the Massey product.
It is possible, however, to obtain a precise realization taking values in $H^2(G_p, \Q_p(1))$
in the present situation, starting from a non-abelian cocycle $a=a_1+a_2$ with values in
$U_2$.
 For this, we make use, on the one hand, of the component $a_2$, which is not a cocycle but satisfies the equation
 $$da_2=-(1/2)(a_1\cup a_1).$$
 That is to say, the cochain $-2a_2$ is one piece of a defining system, already included in
 the cocycle $a$.
 On the other, the equation for the other component $b$ of a defining system looks
 like
 $$db=c^p\cup a_1.$$
 At this point, the assumption on the elliptic curve comes into play implying that  $a_1$ is the localization at $p$ of
 a global cocycle $$a_1^{glob}$$ with the property that the localizations
 $a_1^{glob, l}$
 at all primes $l\neq p$ are trivial. Since $c^p$ is also the localization of
 the global $p$-adic log cyclotomic character $c$, we get an equation
 $$db=c\cup a_1^{glob}$$
 to which we may now seek a global solution, i.e., a  cochain $b$ on a suitable global Galois group. The main point
 then is
 the deep result of Kolyvagin \cite{kolyvagin} that can be used to deduce the existence of
  a global solution $$b^{glob}.$$
 Having solved the equation globally, we again localize to a cochain
 $b^{glob,p}$ on the local Galois group at $p$.
 It is then a simple exercise to see that the Massey product
 $$b^{glob,p}\cup a_1+c^p\cup (-2a_2)$$ obtained thereby
 is independent of the choice of $b^{glob}$, giving us a well-defined function
 $$\psi^p:H^1_f(G_p,U_2) \ra H^2(G_p, \Q_p(1))\simeq \Q_p.$$
 This function is in fact non-zero and algebraic with respect to the structure of $H^1_f(G_p,U_2)$
 as a $\Q_p$-variety. The main theorem then says:
 \begin{thm}
 $$\psi^p\circ j^{et}_{2,\loc}$$
 vanishes on the global points $\cX(\Z)\subset \cX(\Z_p)$.
 \end{thm}
 This result suffices to show the finiteness of integral points. However, the matter
  of real interest is an explicit computation of the composition
 $\psi^p\circ j^{et}_{2,loc}$. At this point, the exponential map
 from the De Rham realization will intervene, and
 the third section provides a   flavor of the formulas one is to expect.
 To get explicit expressions
 we choose a global regular differential form $\a$ on $E$ and a differential $\b$ of
 the second kind with a pole only at $e\in E$ and with the property that
  $[-1]^*(\b)=-\b$. A tangential base-point $b\in T_eE$
 for the fundamental group of $X$
 then determines analytic functions
 $$\log_{\a}(z)=\int_b^z\a, \ \ \log_{\b}(z)=\int_b^z\b, \ \ \ D_2(z)=\int_b^z\a \b,$$
 on $\cX(\Z_p)$ via (iterated) Coleman integration.
 \begin{cor}
 Suppose there is a  point $y$ of infinite order in $\cE(\Z)$. Then
 $$\cE(\Z)\subset \cE(\Z_p)$$
 is in the zero set of
 $$(\log_{\a}(y))^2(D_2(z)-(1/2)\log_{\a}(z)\log_{\b}(z))-
 (\log_{\a}(z))^2(D_2(y)-(1/2)\log_{\a}(y)\log_{\b}(y)).$$
 \end{cor}
We obtain thereby a rather harmonious  constraint on the locus of global integral points, albeit in an absurdly
 special situation.
 In fact, the theorem itself implies that the function of $z$
 $$(D_2(z)-(1/2)\log_{\a}(z)\log_{\b}(z))-
\frac{(D_2(y)-(1/2)\log_{\a}(y)\log_{\b}(y))}{\log_{\a}(y)^2 }(\log_{\a}(z))^2$$
 is independent of the choice of $y$. However, in its present formulation,
 it requires us to have in hand one integral point before commencing the search
 for others.

 Perhaps naively, the author has believed for some time that a satisfactory
 description of the set of global points is possible even for general  hyperbolic curves, compact or affine,
 by way of a non-abelian Poitou-Tate duality of sorts, coupled to an non-abelian explicit reciprocity law
 (cf. \cite{kato}).
 As yet, a plausible formulation of such a duality is unclear, and more so the prospect of
  applications to Diophantine problems.
 What is described in the following sections  is a faint projection of
 the phenomenon whose general nature remains elusive, a projection made possible through the most
 stringent assumptions that are compatible still with statements
 that are not entirely trivial. For the tentative nature of this exposition then,
 even more apologies are in order.

 \section{Preliminary remarks}
 Within the strictures of the present framework, it will be sufficient
 to work with the Selmer variety associated to $U_2=U^3\bsl U$, the first non-abelian level
 of the $\Q_p$-pro-unipotent fundamental group $$U:=\piu$$
 of $$\bX=\times_{\Spec(\Q)}\Spec(\bQ)$$ with a rational tangential base-point $b$
 pointing out of the origin of $E$ \cite{deligne}. Here, the superscript refers to the descending
 central series (in the sense of pro-algebraic groups)
 $$U^1=U, \ \ \ U^{n+1}=[U,U^n]$$
 of $U$ while the subscript denotes the corresponding quotients
 $$U_n=U^{n+1}\bsl U.$$
 In particular, there are canonical isomorphisms
 $$U_1\simeq H_1(\bX,\Q_p)\simeq H_1(\bE,\Q_p)\simeq T_pE\otimes \Q_p$$
 and an exact sequence
 $$0\ra U^3\bsl U^2\ra U_2 \ra U_1\ra 0.$$
 The commutator map induces an anti-symmetric pairing
 $$U_1\otimes U_1 \ra U^3\bsl U^2,$$
 which therefore leads to an isomorphism
 $$U^3\bsl U^2\simeq \wedge^2 H_1(\bar{E}, \Q_p)\simeq \Q_p(1)$$
 as representations for $G=\Gal(\bQ/\Q)$.
 The logarithm map
 $$\log:U \ra L:=Lie U$$
 is an isomorphism of schemes allowing us to identify
 $U_2$ with $L_2=L^2\bsl L$, which, in turn, fits into an exact sequence
 $$0\ra L^3\bsl L^2 \ra L_2 \ra L_1 \ra 0.$$
 If we choose any elements $A$ and $B$  of $L_2$ lifting a basis of
 $L_1$, then $C:=[A,B]$ is a basis for $L^3\bsl L^2$. Using the Campbell-Hausdorff formula,
 we can  express the multiplication in $U_2$, transferred over to $L$ via the logarithm,
 as
 $$(aA+bB+cC)*(a'A+b'B+c'C)=(a+a')A+(b+b')B+(c+c'+(1/2)(ab'-ba'))C.$$
 Given such a choice of $A$ and $B$, we will also denote an element of $L_2$ as
 $l=l_1+l_2$ where $l_1$ is a linear combination of $A$ and $B$ while
 $l_2$ is multiple of $C$. In this notation, the group law becomes
 $$(l_1+l_2)*(l'_1+l'_2)=l''_1+l''_2,$$
 where $l''_1=l_1+l_1'$ and $l''_2=l_2+l'_2+(1/2)[l_1,l'_1]$.
 We simplify notation a bit and put $Z:=L^3\bsl L^2$.
 The Lie bracket
 $$[\cdot, \cdot]:L_2\otimes L_2 \ra Z$$
 factors to a bilinear map
 $$L_1\otimes L_1 \ra Z,$$
 which we denote also by a bracket.

 \begin{lem} Let $p$ be an odd prime. There is a $G$-equivariant vector space splitting
 $$s:L_1\hra L_2$$
 of the exact sequence
 $$0\ra Z \ra L_2 \stackrel{f}{\ra} L_1\ra 0.$$
 \end{lem}
 {\em Proof.}
 Denote by $i$ the involution on $E$ that send $x$ to $-x$ for the group law.
 The origin is fixed and $i$ induces the antipode map
 on the tangent space $T:=T_e(E)$. In particular, we get an isomorphism
 $$i_*:\piu\simeq \mpiu.$$
 Consider  the $\pih$-torsor of paths (for the pro-finite fundamental group)
 $\pibmb$ from $b$ to $-b$. By definition,
 we have
 $$\pibmb:=\Isom(F_b,F_{-b}).$$
 Recall briefly the definition
 of $F_v$  for  $v\in T^0=T\setminus \{0\}$ (\cite{deligne},section 15).
 $F_v$ associates to any cover $Y \ra \bX$, the fiber over $v$ of
 the corresponding cover (`the principal part')
 $$Pr(Y) \ra \bar{T}^0=T^0\otimes \bQ,$$
 of $\bT^0$. Now, choose an isomorphism
 $z:T\simeq \A^1$ that takes $b$ to $1$. Then we get an isomorphism
 $(T^0,b)\simeq (\Gm,1)$ and the pro-$p$ universal covering
 of $\bT^0$ is the pull-back of the tower $$(\cdot)^{p^n}:\bar{\mathbb{G}}_m \ra \bar{\mathbb{G}}_m.$$
 But then, the inverse image of $-1$,
 that is $z^{-1}(-1)$, in each level of the tower gives a compatible $G$-invariant sequence
 of elements and  trivializes the torsor
 $\hpi_1(\bT^0;b, -b)$. Hence, its image in $\pibmb$ will also be a trivialization.
 We then take the unipotent image of this trivialization to get an isomorphism
 $t:\mpiu \simeq \piu$. Note that in the abelianization, we have the canonical isomorphisms
 $$\mpiu^{ab} \simeq H_1(\bE,\Q_p) \simeq \piu^{ab}.$$
 But then $t$, being given by  a path, induces the identity on
 $H_1(\bE,\Q_p)$. Therefore, we have constructed an isomorphism
 $I:=t\circ i_*:U_2\simeq U_2$
 that lifts the map $x\mapsto -x$ on $U_1$.
 This gives a corresponding Lie algebra isomorphism
 $$L_2\simeq L_2,$$ which we will denote by the same letter $I$.
 Since $Z\subset L_2$ is generated by the bracket of two basis elements from
 $L_1$, we see that $I$ restricts to the identity on $Z$.

 Now we define the splitting by putting
 $$s(x):=(1/2)(x'-I(x'))$$
 for any lift $x'$ of $x$ to $L_2$. Since another lift will differ from $x'$ by an element
 of $Z$ on which $I$ acts as the identity, the formula is independent of the lift.
 We can then use this independence to check that the map is linear. If $x'$ and $y'$ are
 lifts of $x$ and $y$, then $\l x'+ \mu y'$ is a lift of $\l x+ \mu y$. So
 $$s(\l x+\mu y)=(1/2)(\l x'+\mu y'-I(\l x'+\mu y'))=\l (1/2)(x'-I(x'))+\mu (1/2)(y'-I(y')).$$
  Similarly, if
 $g\in G$, then $g(x')$ will be a  lift of  $g(x)$, so that
 $$s(g(x))=(1/2)(g(x')-I(g(x')))=(1/2)(g(x')-g(I(x'))=g(s(x)).$$
 $\Box$

 We will use this splitting to write
 $$L_2=L_1\oplus Z$$
 as a $G$-representation, so that  an arbitrary $l\in L_2$  can be decomposed into
 $l=l_1+l_2$ as described above, except independently of a specific basis of $L_1$.
  Using the identification via the log map, we will  abuse notation
 a bit and write $u=u_1+u_2$ also for an element of $U_2$.

For any $\l \in \Q_p$, we then get a Lie algebra homomorphism
$$m(\l): L_2\simeq L_2$$
by defining
$$m(\l) l=\l l_1+\l^2 l_2,$$
which is furthermore compatible with the $G$-action. Thus, $m(\l)$ can also be
viewed as a $G$-homomorphism of $U_2$. We note that the extra notation is used
for the moment to distinguish this action of the multiplicative monoid $\Q_p$
from the original scalar multiplication.

 We recall the continuous group cohomology \cite{kim1}
 $$H^1(G,U_2):=U_2\bsl Z^1(G,U_2).$$
 The set $Z^1(G,U_2)$ of continuous 1-cocycles of $G$ with values in $U_2$
 consist of continuous maps
 $$a:G \ra U_2$$
 such that
 $$a(gh)=a(g)ga(h).$$
 Using the identification of
 $U_2$ with $L_2=L_1\oplus Z$ just discussed, we will write
 such a map also as
 $$a=a_1+a_2$$
 with $a_1$ taking values in $L_1$ and $a_2$ values in $Z$.
 The cocycle condition is then given by
 $$a_1(gh)+a_2(gh)=(a_1(g)+a_2(g))*(ga_1(h)+ga_2(h))$$
 $$=a_1(g)+ga_1(h)+a_2(g)+ga_2(h)+(1/2)[a_1(g),ga_1(h)].$$
 In fact, given two cochains $c,c'$ with values in $L_1$, we define
 $$(c\cup c')(g,h)=[c(g),gc'(h)],$$ a cochain with values in $Z$.
 The cocycle condition spelled out above then says that
 \bq
 (1) $a_1$ is a cocycle with values in $L_1$;

 (2) $a_2$ is not a cocyle in general, but satisfies
 $$ga_2(h)-a_2(gh)+a_2(h)=-(1/2)[a_1(g),ga_1(h)],$$
 \eq
 or, recognizing on the left hand side the differential of the cochain $a_2$,
 $$da_2=-(1/2)(a_1\cup a_1).$$
 This can be viewed as a Galois-theoretic {\em  Maurer-Cartan equation} \cite{GM}.
 So even when we have split the Galois action at this first-nonabelian level, the non-abelian group
 structure imposes a twist on cohomology. The cohomology set is then defined by
 taking a quotient under the action of $U_2$:
 $$(u\cdot a)(g)=ua(g)g(u^{-1}).$$
 This discussion carries over verbatim to various local Galois group $G_l=\Gal(\bQ_l/\Q_l)$ as $l$ runs over the
 primes of $\Q$, which all act on $U$ via the inclusion $G_l \ra G$ induced by an
 inclusion $\bQ \hra \bQ_l$, and leading to the cohomology sets $H^1(G_l,U_2)$.

Now we choose $p$ to be an odd place of good  reduction for $E$ and
 denote by $T$ the set $S\cup \{p\}$, where $S$ is  a finite set of places that contains the infinite place
 and the places of bad reduction
 for $E$. We let $G_T=\Gal(\Q_T/\Q)$ be the Galois group of the maximal extension of
 $\Q$ unramified outside $T$.
 In previous work  \cite{CK, kim1, kim2, kim3, kim4, KT}, we have made use of the Selmer variety
 $$H^1_f(G,U_2)\subset H^1(G_T,U_2)\subset  H^1(G,U_2).$$
 By definition,
 $H^1_f(G,U_2)$  consists of cohomology classes that
 are unramified outside $T$ and crystalline at $p$.
 Associating  to a point $x\in X(\Q)$ the torsor
 of paths
 $$\Pux$$
 defines a map
 $$X(\Q) \ra H^1(G,U_2)$$
 that takes $\cX(\Z_S)\subset X(\Q)$ to
 $H^1_f(G,U)$.
 In fact, let $l\neq p$, and $G_l:=\Gal(\bQ_l/\Q_l)$.
 Then on the globally integral points, the map
 $$\cX(\Z) \ra H^1_f(G,U_2) \ra H^1(G_l,U_2)$$
 is trivial. This was essentially shown in \cite{KT}. What is shown there
 is that the image factors through the map
 $$H^1(G_l/I_l,U^{I_l}_2)\ra H^1(G_l,U_2),$$
 where $I_l\subset G_l$ is the inertia subgroup. But for both
 $U_1$ and $Z$, the Frobenius element in $G_l/I_l$ acts with strictly negative weights,
 and hence, $$H^1(G_l/I_l, Z^I)=H^1(G_l/I_l, U_1^I)=0,$$
 from which we deduce that
 $$H^1(G_l/I_l,U^{I_l}_2)=0.$$
Denote by
$$H^1_{f,\Z}(G,U_2)\subset H^1_f(G,U_2)$$
the {\em fine Selmer variety}, defined to be the  intersection of the kernels of the localization maps
$$H^1_f(G,U_2) \ra H^1(G_l,U_2)$$
for all $l\neq p$. Then the previous paragraph says that the image of
$$\cX(\Z) \ra H^1_f(G,U_2)$$
lies inside
$$H^1_{f,\Z}(G,U_2).$$

 The action of $\Q_p$ discussed above induces an
 action on $H^1_{(\cdot)}(\cdot, U_2)$ for any of the groups under discussion, which we will denote
 simply as left multiplication (since here, the danger of confusion
 with the original scalar multiplication does not arise). At the level of cocycles,
 $$\l a= \l a_1+\l^2 a_2.$$
 \section{Construction}
 \smallskip

 {\em (i) Global construction}
 \smallskip

Let $\chi:G_T \ra \Z_p^*$ be the $p$-adic cyclotomic character
and $c:=\log \chi: G_T
 \ra \Q_p,$ regarded as an element of
$H^1(G_T, \Q_p)$. For any point $x\in \cX(\Z_S)$, let $a(x)$ be a cocycle representing the
class of $\pi_1^{u,\Q_p}(\bX;b,x)$ in $H^1(G_T,U_2)$.
Write
$$a(x)=a_1(x)_1+a_1(x)$$
as indicated in the previous section.
Then $c\cup a_1(x)$ represents a cohomology class in
$H^2(G_{T}, L_1).$
But if we consider the localization map
$$0\ra \Sha^2_{T}(L_1)\ra H^2(G_{T},L_1) \ra \oplus_{v\in {T}}H^2(G_v,L_1),$$
we see that
$$H^2(G_v,L_1)\simeq H^0(G_v,L_1)^*=0$$
for all $v$.
On the other hand, the kernel $\Sha^2_{T}(L_1)$  is dual to the kernel $\Sha_{T}^1(L_1)$ of
the $H^1$-localization
$$0\ra \Sha_{T}^1(L_1)\ra H^1(G_T,L_1) \ra \oplus_{v\in T}H^1(G_v,L_1),$$
which then must  lie inside the $\Q_p$-Selmer group
$$H^1_{f,0}(G,L_1):=(\invlim_n Sel(E)[p^n])\otimes \Q_p.$$ However, because of our assumption $\ord_{s=1} L(E,s)=1$, we know that
the Tate-Shafarevich group of $E$ is finite \cite{kolyvagin}, and that
$$E(\Q)\otimes \Q_p\simeq H^1_{f,0}(G,L_1).$$
Hence, there is an injection
$$H^1_{f,0}(G,L_1)\hra H^1(G_p,L_1),$$
from which we deduce that $\Sha^1_T(L_1)=0$.
We conclude that $H^2(G_T,L_1)$ vanishes, allowing us to find a cochain $b_x:G_T \ra L_1$ such that
$$db_x=c\cup a_1(x).$$
Define a two-cochain
$$\phi_x:G_T\times G_T\ra Z$$
by putting
$$\phi_x=b_x\cup a_1(x)-2c\cup a_2(x).$$

\begin{lem}
$\phi_x$ is a cocycle.
\end{lem}
{\em Proof.} Since $a(x)_1$ and $c$ are cocycles,
we have
$$d\phi_x=db_x\cup a(x)_1+2c\cup da(x)_2=c\cup a(x)_1\cup a(x)_1-c\cup a(x)_1\cup a(x)_1=0.$$
$\Box$

Our construction of $\phi_x$ depends on the auxiliary cochain $b_x$, and hence, gives a
class
$$[\phi_x]\in H^2(G_T, Z)/[H^1(G_T,L_1)\cup a_1(x)].$$
\begin{lem}
The class $[\phi_x]$ is independent of the choice of cocycle $a(x)$.
\end{lem}
{\em Proof.} Obviously, the subspace
$H^1(G_T,L_1)\cup a_1(x)$ depends only on the class of $a_1(x)$, and hence,
on the class of $a(x)$.
Now we examine the action of $U_2$. To reduce clutter, we will
temporarily suppress $x$ from the notation for the cochains.
Write $u=u_1+u_2$.
Then
$$ua(g)g(u^{-1})=(u_1+u_2)*(a_1(g)+a_2(g))*(-gu_1-gu_2)$$$$=(u_1+a_1(g)+u_2+a_2(g)+(1/2)[u_1,a_1(g)])*(-gu_1-gu_2)$$
$$=a_1(g)+u_1-gu_1+a_2(g)+u_2-gu_2+(1/2)[u_1,a_1(g)]-(1/2)[a_1(g),gu_1]-(1/2)[u_1,gu_1].$$
The $L_1$-component of this expression is
$a_1-du_1$, where we view the element $u_1\in L_1$ as a zero-cochain.
Thus, the previous choice of $b_x$ can be changed to
$b_x+c\cup u_1$, since
$$d(c\cup u_1)=dc\cup u_1-c\cup du_1=-cdu_1.$$
The resulting two-cocycle changes to
$$(b_x+c\cup u_1)\cup (a_1-du_1)-2c\cup a_2+2c\cup du_2-c\cup r+c\cup s+c\cup t$$
$$=\phi_x+[c\cup u_1\cup a_1- b_x\cup du_1-c\cup u_1\cup du_1+2c\cup du_2-c\cup r+c\cup s+c\cup t],$$
where $r,s,t$ are the functions $G_T\ra Z$ defined by
$$r(g)=[u_1,a_1(g)], \ \ s(g)=[a_1(g),gu_1],\ \ t(g)=[u_1,gu_1].$$
Within the discrepancy, the term $2c\cup du_2=-2d(c\cup u)$ is clearly a co-boundary. Also,
we see that $t=u_1\cup du_1$, ridding us of two  terms.
We have
$$c\cup u_1\cup a_1(g,h)=[c\cup u_1(g),ga_1(h)]=[c(g)gu_1, ga_1(h)],$$
while
$$c\cup r (g,h)=c(g)gr(h)=c(g)g[u_1,a_1(h)]=c(g)[gu_1,ga_1(h)],$$
causing two more terms to cancel each other.
Finally, it remains to analyze the difference
$c\cup s-b_x\cup du_1$.
But
$$d(b_x\cup u_1)=db_x\cup u_1-b_x\cup du_1,$$
so that up to a co-boundary, we can replace
$b_x\cup du_1$ by $db_x\cup u_1=c\cup a_1\cup u_1$.
We can then compute the value
$$c\cup a_1\cup u_1(g,h)=[c\cup a_1(g,h), ghu_1]=[c(g)ga_1(h), ghu_1],$$
which is verified to be equal to
$$c\cup s (g,h)=c(g)gs(h)=c(g)g[a_1(h),hu_1]=c(g)[ga_1(h),ghu_1].$$
Therefore, the difference
$c\cup s-b_x\cup du_1$ is a coboundary.
$\Box$

Given a global cohomology class or cochain $s$, we will denote by $s^l$ its localization
at the prime $l$, that is, its restriction to $G_l$.

\begin{lem}
The subspace
$H^1(G_T,L_1)\cup a_1(x)$ of $H^2(G_T,Z)$ is zero.
\end{lem}
{\em Proof.}
Because $a_1(x)$ is the class of a point and we are taking
$\Q_p$-coefficients, $a_1(x)^l=0$ for all $l\neq p$. Thus, for any
$r\in H^1(G_T,L_1)$, we have $(r\cup a_1(x))^l=0$ for
all $l\neq p$. This is of course also true for the archimedean component since
$p$ is odd. Thus, the only possible component of $r\cup a_1(x)$ that survives
is at $p$, which then must be zero since
$$\sum_v(r\cup a_1(x))_v=0.$$
Because we also have an injection
$$0\ra H^2(G_T,\Q_p(1)) \hra \oplus_v H^2(G_v, \Q_p(1)),$$
we see that $r\cup a_1(x)=0$.
$\Box$

We actually have therefore a well-defined class
$$[\phi_x]\in H^2(G_T,Z).$$

\begin{lem}
Let $l\neq p$. Then
$$[\phi_x]^l\in H^2(G_l, Z)$$
can be computed locally in the following sense:
Choose any local representative $t(x)$ for the class
$[a(x)]^l$, and let
$s:G_l\ra L_1$ be any local cochain such that
$ds=c^l\cup t_1(x)$. Then
$$[\phi_x]^l=[s\cup t_1(x)-2c^l\cup t_2(x)].$$
\end{lem}
{\em Proof.}
The class modulo
$$H^1(G_l,L_1)\cup t_1(x)$$
will clearly be independent of the choice of $t $ and $s$.
Since $[t_1(x)]=[a^l_1(x)]$, we have
$$H^1(G_l,L_1)\cup t_1(x)=H^1(G_l,L_1)\cup a_1(x).$$
But, as we pointed out above, $a^l_1(x)$ is the trivial class.
Therefore, the class in $H^2(G_l, Z)$ is independent of the choices.
In particular, local choices $s$ and $t$ will give the same class
as the localization of the global choices $a$ and $b$.
$\Box$
\medskip

{\em (ii) Local construction}
\smallskip

We will now make use of a point $y\in E(\Q)$ of infinite order.
 Given any class
$s\in H^1_f(G_p,U_2)$, its component
$s_1\in H^1_f(G_p, U_1)$ is a $\Q_p$ multiple of
$a^p_1(y)$:
$$s_1=\l (s) a^p_1(y),$$
for some $\l (s)\in \Q_p$. In particular,
$s_1=x_1^p$ for some
cocycle $x_1:G_T\ra L_1$ such that
$[x_1]^l=0$ for all $l\neq p$.
By the theorem of Kolyvagin cited earlier, the equation
$$db=c\cup (\l(s) x_1)$$
has a solution $b^{glob}:G_T\ra L_1$.
Thus, we get a class
$$\psi^p(s):=[b^{glob,p}\cup s_1-2c^p\cup s_2]\in H^2(G_p, Z)/[\loc_p(H^1(G_T,L_1))\cup s_1]$$
since two  choices of $b^{glob}$ will differ by an element of $H^1(G_T,L_1)$.
But
$$\loc_p(H^1(G_T,L_1))\cup s_1=\loc_p(H^1(G_T,L_1))\cup a^p_1(y)=\loc_p((H^1(G_T,L_1))\cup a_1(y))=0$$
by Lemma 2.3. Therefore, we have a well-defined class
$$\psi^p(s)\in H^2(G_p,Z).$$

The following lemma is straightforward from the definitions:
\begin{lem}
Let $x\in \cX(\Z_S)$. Then
$$\psi^p(a^p(x))=[\phi_x]^p.$$
\end{lem}
\medskip

Now we can give the
\medskip

{\em Proof of theorem 0.1}

If $x\in \cX(\Z)$, then we know that
$[a^l(x)]=0$ for all $l\neq p$. By lemma 2.4, this implies that
$[\phi_x]^l=0 $ for all $l\neq p$, and hence,
$$\psi^p(a^p(x))=[\phi_x]^p=0.$$
$\Box$
\medskip

Any class $s\in H^1_f(G_p,U_2)$  lies over the same point in
$H^1_f(G_p,U_1)$ as $\l(s) a^p(y)$ \footnote{Take care that this multiplication now refers to the
$\Q_p$-action discussed in the previous section.}. By the exact sequence \cite{KT}
$$0\ra H^1_f(G_p, Z)\ra H^1_f(G_p,U_2) \ra H^1_f(G,U_1),$$
the two classes then differ by the action of an element of
$H^1_f(G_p, Z)$ which we denote by
$$\l(s)  a^p (y)-s \in H^1_f(G_p, Z).$$
Using the point $y$, we get the following alternative description of
the function $\psi$:
\begin{lem}
$$\psi^p(s)=\l(s)^2\psi^p(y)+2(c^p\cup (\l(s)a^p(y)-s)).$$
\end{lem}

{\em Proof.}
Let $b:G_T\ra \L_1$ be a solution of
$$db=c\cup a_1(y).$$
Then
$$d(\l(s) b)=c\cup (\l(s) a_1(y))$$
and
$$(\l(s) b)^p\cup s_1=(\l(s) b)^p\cup (\l(s) a^p_1(y))=\l(s)^2 b^p\cup a^p_1(y).$$
Therefore,
$$\psi^p(s)=\l(s)^2 b^p\cup a^p_1(y)-2c^p\cup s_2$$
$$=\l(s)^2 (b^p\cup a^p_1(y)-2c^p\cup a^p_2(y))+
2(\l(s)^2 c^p\cup a^p_2(y)-c^p\cup s_2)$$
$$=\l(s)^2 \psi^p(y)+2(c^p\cup (\l(s)a^p(y)-s)).$$
$\Box$
\medskip

Fix now the isomorphism
$Z\simeq \Q_p(1)$ induced by the Weil pairing $<\cdot, \cdot>$, that is,
that takes
$[x,y]$ to $<x,y>$, which then induces an isomorphism
$$T: H^2(G_p, Z)\simeq \Q_p$$
and gives us a $\Q_p$-valued function
$$T\circ \psi^p:H^1_f(G_p,U_2) \ra \Q_p.$$
We will sometimes suppress $T$ from the notation and simply regard $\psi^p$
as taking values in $\Q_p$.
From the definition, we see that for any $\l \in \Q_p$,
$$\psi^p_y(\l a^p(y))=\l^2 T(\phi^p_y),$$
while for $r\in H^1_f(G_p, Z)$, we have
$$\psi^p_y(r)=-T(c^p\cup r).$$
Thus, when we take $z\in \Z_p^*$ with a cohomology class
$$k(z)\in H^1_f(G_p, \Q_p(1))$$
coming from Kummer theory that we  identify with a class in
$H^1_f(G_p,Z)$, then
$$\psi^p_y(z)=-\log \chi(Rec_p(z)),$$
where $Rec_p$ is the local reciprocity map. Since
$\chi (Rec_p(z))=z$, we get
$$\psi^p_y(k(z))=-\log z\in \Q_p.$$
In particular, the map is not identically zero. In fact, it is far from trivial
on any  fiber of
$$H^1_f(G_p,U_2)\ra H^1_f(G_p,U_1).$$

With respect to the structure of $H^1_f(G_p, U_2)$ as an algebraic variety,
$\psi^p$ is in fact a non-zero algebraic function, as we see in a straightforward way by defining it for points in arbitrary
$\Q_p$-algebras (as in \cite{kim1}). Since there is
a Coleman map
$$j^{et}_{2,loc}:\cX(\Z_p)\ra H^1_f(G_p,U_2)$$
with Zariski dense image for each residue disk \cite{kim2},  Theorem 0.1 yields the finiteness of
$\cX(\Z)$. As mentioned in the introduction, the obvious task of importance is to compute the function
$$\psi^p\circ j^{et}_{2, loc}$$
on $\cX(\Z_p)$, whose zero set is guaranteed to capture the global integral points.

\section{Preliminary formulas}

There is a commutative diagram \cite{kim2}
$$\bd \cX(\Z_p) & \rTo^{j^{et}_{2,loc}}& H^1_f(G_,U_2) \\
 & \rdTo^{j^{dr/cr}_2} & \dTo^{\simeq} \\
 & & U^{DR}_2/F^0 \ed$$
 bringing  the De Rham fundamental group $U^{DR}=\pi_1^{DR}(X_{\Q_p},b)$ and its quotient $U^{DR}_2=U^{DR}/U^3$
 into our consideration.
 The isomorphism
 $$H^1_f(G_,U_2) \simeq U^{DR}_2/F^0$$
 is a non-abelian analogue of the Bloch-Kato log map.
 There is actually a larger commutative diagram
$$\bd \cX(\Z_p) & \rTo^{j^{et}_{loc}}& H^1_f(G_,U) \\
 & \rdTo^{j^{dr/cr}} & \dTo^{\simeq} \\
 & & U^{DR}/F^0 \ed$$
 out of which the level-two version is obtained by composing with the projection
 $$U\ra U_2.$$
 The map
 $j^{et}_{loc}$ is just a local version of the map recalled in the previous section, while
 $j^{dr/cr}$ associates to each point $x\in \cX(\Z_p)$, the $U^{DR}$-torsor
 $\Pdrx$, and hence, corresponds to a point $j^{dr/cr}(x)\in U^{DR}/F^0$ \cite{kim2}.
 The point is described explicitly as follows.
 One chooses an element $p^H\in F^0\Pdrx$. On the other hand, there is a unique
 Frobenius invariant element $p^{cr}\in \Pdrx$. Then $j^{dr/cr}(x)$ the coset of the element
 $u$ such that $p^{cr}u=p^H$.
 In \cite{kim2} we gave a description of a universal pro-unipotent bundle with connection
 on $X_{\Q_p}$. Let $\a$ be an invariant differential 1-form on $E$ and let $\b$
 be a differential of the second kind with a pole only at $e$ such that $[-1]^*(\a)=-\a$ and $[-1]^*(\b)=-\b$. (Of course the
 first condition is automatic.)
  Let
 $$R:=\Q_p<<A,B>>=\invlim \Q_p<A,B>/I^n,$$
 where $\Q_p<A,B>$ is the free-noncommutative $\Q_p$-algebra on the letters $A,B$,
 and $I\subset \Q_p<A,B>$ is the augmentation ideal. Thus, $\Q_p<A,B>$ is spanned
 by words $w$ in $A$ and $B$, while
 the elements of $R$ are infinite formal linear combinations
 $$\sum c_w w$$
 in such words with coefficients $c_w\in \Q_p$.
 Using this we can construct
 the free $\O_{X_{\Q_p}}$-module
 $\cR=\O_{X_{\Q_p}} \otimes R$ together with the connection
 $$\nabla f =df-(A\a+B\b) f,$$
 for an element $f\in \cR$. If we choose the element $1\in \cR_b=R$ as the initial condition,
 then the element $p^{cr}(1)\in \cR_x$ corresponding to it is given by
 $$
 G(x)=\sum_w \int_b^x a_w  w,$$
 where $a_w$ is the symbol
 $$\a^{n_1}\b^{m_1}\cdots \a^{n_k}\b^{m_k}$$
 if $w$ is the word
 $$A^{n_1}B^{m_1}\cdots A^{n_k}B^{m_k},$$
 and the integral symbol corresponds to iterated Coleman integration (\cite{furusho}, \cite{BF} Prop. 4.5). This is normalized
 by the convention
 $$d(\int \a a_w)=(\int a_w )\a,\ \ \ d(\int \b a_w)= (\int a_w )\b.$$
 To recall the detailed description, we need to choose a local coordinate $z$ at $e\in E$ such that
 $d/dz=b$. Furthermore, fix Iwasawa's branch of the $p$-adic log. This determines a ring
 $A^{col}(]e[)(\log z)$ of logarithmic
 Coleman functions in the residue disk $]e[$ of the origin of $E$, and
 there is a unique element (see op. cit.)
 $G^b\in A^{col}(]e[)(\log z)\otimes R$ characterized by the property
 $\nabla G=0$, together with the initial condition specifying that when we write
 $$G^b=G^b_0+G^b_1\log z+G^b_2(\log z)^2+\cdots,$$
 we have $G^b_0(0)=1$.
 Then analytic continuation along the Frobenius produces an element
 $$G^x\in A^{col}(]x[)\otimes R$$ in the residue disk $]x[$ compatible with $G^b$, and
 $$G(x):=G^x(x).$$
 There is a co-multiplication
 $$\D:R\ra R\otimes R$$
 determined by
 $$\D(A)=A\otimes 1+1\otimes A, \ \ \D(B)=B\otimes 1+1\otimes B,$$
 with respect to which $G(x)$ is group-like, i.e., satisfies
 $$\D(G(x))=G(x)\otimes G(x).$$
 Thus, $G(x)$ corresponds to a $\Q_p$ point of
 $\Pdrx=\Spec(R^*)$, where
 $$R^*:=\dirlim \Hom(R/I^n,\Q_p).$$
 The structure of $R^*$ can also be elucidated as the
 free $\Q_p$-vector space generated by the functions
 $\a_w$ such that $\a_w(w')=\d_{w w'}$.
 Another description of
 $R^*$ is as the $H^0$ of the circular reduced  bar construction (\cite{HZ}, section 3)
 $$H^0(B( TW(G(\Om^*_{E_{\Q_p}}(\log e))), \Q_p)),$$
 where
 $TW(G(\Om^*_{E_{\Q_p}}(\log e)))$ is the Thom-Whitney construction $TW(\cdot)$
  on the Godement resolution $G(\cdot)$ of $\Om_{E_{\Q_p}}(\log e)$ \cite{KH, navarro}, the sheaf of differential forms on
 $E_{\Q_p}$ with log poles along $e$, and the $\Q_p$ is the $TW(G(\Om^*_{E_{\Q_p}}(\log e)))$-bi-module obtained by
 evaluation at $b$ and $x$.
 We have a multiplicative quasi-isomorphism of sheaves of commutative differential graded algebras \cite{deligne2}
 $$\Om^*_{E_{\Q_p}}(\log e) \ra j_*\Om^*_{X_{\Q_p}},$$
 where $j:X\hra E$ is the inclusion. The Hodge filtration $F$ on the left is compatible with
 the pole-filtration $P$ on the right, where
 $$P^0j_*\Om^*_{X_{\Q_p}}=[\O_{E_{\Q_p}}(e) \ra \Om_{E_{\Q_p}}(2e)],$$
 and
 $$P^1j_*\Om^*_{X_{\Q_p}}=[0 \ra \Om^1_{E_{\Q_p}}(e)].$$
 Thus, we have filtered quasi-isomorphisms
 $$TW(G(\Om^*_{E_{\Q_p}}(\log e)))\simeq TW(G(j_*\Om^*_{X_{\Q_p}})).$$
 But the latter admits a filtered quasi-isomorphism
 $$\Om^*_{X_{\Q_p}}(X_{\Q_p}) \ra TW(G(j_*\Om^*_{X_{\Q_p}})),$$
 where the pole filtration on the left hand side has simply
 $$P^1\Om^*_{X_{\Q_p}}(X_{\Q_p})=[0\ra \Om^1_{E_{\Q_p}}(E_{\Q_p})].$$
 Therefore, $F^1$ of the Hodge filtration on
 $R^*$ is generated by the $\a_w$ such that $w$ contains at least one $A$.
 Hence, the Hodge filtration on $R$ has $F^0$ generated by the
 elements $B^n$.
 Therefore, with respect to the basis $\{A,A^2,AB, BA\}$
 of
 $U^{DR}_2/F^0$, we can express $j^{dr/cr}_2$ as
 $$j^{dr/cr}_2(x)=1+\int_b^x \a A+\int_b^x \a^2 A^2+\int_b^x \a \b AB +\int_b^x \b \a BA.$$
 This map as well is conveniently expressed  in terms of the logarithm
 $\log U^{DR}\ra L^{DR}$as
 $$\log j^{dr/cr}(x)= \int_b^x \a A+(\int_b^x \a \b -(1/2)(\int_b^x\a)(\int_b^x \b))[A,B].$$
Introducing the notation
 $$\log_{\a}(x)=\int_b^x\a, \ \ \log_{\b}(x)=\int_b^x\b, \ \ D_2(x)=\int_b^x\a \b,$$
 we can also write this as
 $$\log j^{DR}(x)=\log_{\a}(x)A+(D_2(x)-(1/2)\log_{\a}(x)\log_{\b}(x))[A,B].$$
 Regarding the $\Q_p$-action, our choice of $\a$ and $\b$ implies that
 the automorphism   of $L^{DR}_2$  induced by the involution of section 1 is simply
 $$A\mapsto -A, \ \ \ B\mapsto -B,$$
 so  that $A$ and $B$ are   basis elements compatible with the grading
 on $L_2$. In particular, the $\Q_p$-action can be described as
 $$m(\l) A= \l A, \ \ \ , m(\l )B=\l B, \ \ \ m(\l)[A,B]=\l^2 [A,B].$$
 Given the class $a(x)\in H^1_f(G_p,U_2)$ of a point $x\in \cX(\Z_p)$, the number
 $\l $ such that $a_1(x)=\l a^p_1(y)$ can  be written as
 $$\l=\log_{\a}(x)/\log_{\a}(y),$$
 since the logarithm is a group homomorphism.
 On the other hand, we have
 $$\log j^{dr/cr}_2(x)=\log_{\a}(x) A+(D_2(x)-(1/2)\log_{\a}(x)\log_{\b}(x))[A,B]$$
 $$\log \l j^{dr/cr}_2(y)=\l \log_{\a}(y) A+ \l^2 (D_2(y)-(1/2)\log_{\a}(y)\log_{\b}(y))[A,B]$$
 $$= \log_{\a}(x) A+ \l^2 (D_2(y)-(1/2)\log_{\a}(y)\log_{\b}(y))[A,B].$$
 So the element in
 $Z^{DR}$ representing the difference is
 $$[(\log_{\a}(x)/\log_{\a}(y))^2(D_2(y)-(1/2)\log_{\a}(y)\log_{\b}(y))
 -(D_2(x)-(1/2)\log_{\a}(x)\log_{\b}(x))][A,B].$$
  Using the (Bloch-Kato) exponential notation
 for the isomorphism
 $$\Exp: U^{DR}_2/F^0\simeq H^1_f(G_p,U_2),$$
 we see then that a formula
 for the function
 $\psi^p\circ j^{et}_{loc}$
 is given by
 $$\psi^p\circ j^{et}_{loc}(x)=(\log_{\a}(x)/\log_{\a}(y))^2T(\psi^p(y))+$$
 $$2[(\log_{\a}(x)/\log_{\a}(y))^2(D_2(y)-(1/2)\log_{\a}(y)\log_{\b}(y))
 -(D_2(x)-(1/2)\log_{\a}(x)\log_{\b}(x))]T(c\cup \Exp([A,B])).$$
 Therefore, obtaining a `concrete' expression for this function reduces to the computation of a single
 $$T(\psi^p(y))$$
 and
 $$T(c\cup \Exp([A,B])).$$
 The latter is found in a rather straightforward manner.
 The key point is that the map
 $$H_1^{et}(\bGmp, \Q_p)\simeq Z\hra \U_2$$
 is induced by the restriction functor
 $$r:Cov(\bX) \ra Cov(\bT^0)$$
 mentioned in section 1.
 Similarly, the map
 $$H_1^{DR}(\Gmp) \simeq Z^{DR}\hra U^{DR}_2$$
 is induced by a restriction functor
 $$r^{DR}:\Un(X) \ra \Un(T^0),$$
 on categories of unipotent bundles with connection.
 The construction of $r^{DR}$ takes a bundle
 $(V,\nabla)$ and associates to it  first
 the canonical extension
 $$(\bar{V}, \bar{\nabla})$$
 which is a log connection
 $$\bar{\nabla}:\bar{V} \ra \Om_E(\log e)\otimes \bar{V},$$
 on $E$, and then its residue
 $$Res(\nabla):V_e\ra V_e,$$
 which is an endomorphism of the fiber $V_e$.
 The value $r^{DR}(V,\nabla)$ is then just the
 trivial bundle $V_e$ on $T^0$ equipped with the connection
 $$d-Nd/dt$$
 for any choice of linear coordinate on $T$ such that $t(0)=0$.
 Let us compute $r^{DR}$ for
 $\cR/I^3$.
 In a formal neighborhood of $e$, we can solve the equation
 $$dv=\b$$
 and make the gauge transformation induced by
 $1-vB$. Then the connection form with respect to this gauge becomes
 $$(1-vB)(\a A+\b B)(1+vB+v^2B^2)+(-dvB)(1+vB+vB^2)$$
 $$=(\a A+\b B-v \a BA-v \b B^2)(1+vB+v^2B^2)-dvB-vdvB^2$$
 $$=\a A+\b B-v \a BA-v \b B^2+v\a AB+v\b B^2-dvB-vdvB^2$$
 $$=\a A+v \a [A,B]-vdvB^2.$$
 modulo $I^3$.
 From this, we see that the residue is
 $$Res(v\a)[A,B]$$
 with $Res(v\a)\in \Z_p^*$.
 This residue
 can be easily identified with the Serre duality pairing $<\a,\b>$.
 In any case, the connection
 $r^{DR}(\cR/I^3)$
 is the bundle
 $\O_{T^0}\otimes R/I^3$
 with connection
 $$d-Res(v\a)[A,B]d/dt.$$
 The universal unipotent connection
 on $T^0$ is
 $$\O_{T^0}\otimes \Q_p[[C]]$$
 with connection form
 $$d-Cd/dt.$$
 So the map
 $$C\mapsto Res(v\a)[A,B]$$
 realizes the map
 $$\pi^{DR}(T^0, b) \ra U^{DR}_2.$$
 In particular, this map fits into the commutative diagram:
 $$\bd
 \pi^{DR}(T^0_{\Q_p}, b) & \rTo &U^{DR}_2/F^0 \\
 \dTo^{\Exp} & & \dTo^{\Exp} \\
 H^1_f(G_p, \pi^{\Q_p,et}_1( \bT^0,b)) & \rTo & H^1_f(G_p, U_2).
 \ed$$
 We analyze the left vertical arrow.
 Fix the linear coordinate $t:T\ra \mathbb{A}^1$ such that $t(0)=0$
 and $t(b)=1$. This induces isomorphisms
 $$\pi^{\Q_p,et}( \bT^0,b))\simeq \Q_p(1)$$
 $$\pi^{DR}(T^0_{\Q_p}, b)\simeq H_1^{DR}(\mathbb{G}_{m,\Q_p}).$$
 The isomorphism
 $$H^1_f(G_p, \Q_p(1))\simeq H_1^{DR}(\mathbb{G}_{m,\Q_p})$$
 takes the class $k(x)$ of $x\in \Z_p^*$ to the class of
 $$\int_1^x dt/t C=\log(x)C.$$
 On the other hand,
 $$T(c\cup k(x))=c(Rec_p(x))=\log \chi(x)=\log x.$$
 Thus, we see that $T(c\cup \Exp(C))=1$.
 Therefore, from the previous commutative diagram, we deduce:
 \begin{prop}
 $$T(c\cup \Exp[A,B])=Res(v\a)^{-1}.$$
 \end{prop}

 The computation of $T(\psi^p(y))$ in general appears to be somewhat difficult. Perhaps some progress is
 possible through the theory of $p$-adic uniformization when $p$ is a split-semi-stable prime, for which
 a generalization of the local Selmer theory will be  necessary. In that case, it will be natural to
 take $y$ the trace of a Heegner point coming from a Shimura curve uniformization, and some
 relationship between the various quantities and $L$-functions should emerge, such as
 $$(\log_{\a}y)^2=(1/2)(d^2/dk^2)L_p(E, k, k/2)_{k=2},$$
 which appears in \cite{BD}.

 One case that is tractable right now is when we already have an integral point $y$ of infinite
 order in hand, because then $T(\psi^p(y))=0$.
 Corollary 0.2 is an immediate consequence.
 \medskip

 {\bf Acknowledgements:} It is a great pleasure to express my
 profound gratitude to John Coates for innumerable mathematical discussions as well as for
 a continuum of moral support surrounding this research.
 Considerable benefit was  derived also from communication with Henri Darmon, Christopher Deninger, and
 Gerd Faltings.

 This paper was completed while the author was a guest of the
 SFB 478 at the University of M\"unster. The excellent academic and cultural
 environment of the city
 was an indispensable aid to the process of writing.

{\footnotesize  Department of Mathematics, University College London,
Gower Street, London, WC1E 6BT, United Kingdom and The Korea Institute
for Advanced Study, Hoegiro 87, Dongdaemun-gu, Seoul 130-722, Korea}

\end{document}